\documentclass[10pt]{amsart}

\usepackage{amssymb,latexsym,amsfonts}

\newcommand{\tensor}{\otimes}

\newcommand{\Hom}{\operatorname{Hom}}

\renewcommand{\c}{\operatorname{c}}

\newcommand{\picard}{\operatorname{Pic}}
\newcommand{\Ext}{\operatorname{Ext}}
\newcommand{\sext}{\operatorname{\mathcal{E}\mathit{xt}}}

\newcommand{\shom}{\operatorname{\mathcal{H}\mathit{om}}}
\renewcommand{\H}{\operatorname{H}}
\newcommand{\h}{\operatorname{h}}

\newcommand{\ses}[3]{0\rightarrow#1\rightarrow#2
   \rightarrow#3\rightarrow0}

\newcommand{\F}{{\mathcal F}}
\newcommand{\G}{{\mathcal G}}

\newcommand{\I}{{\mathcal I}}

\renewcommand{\O}{{\mathcal O}}
\renewcommand{\P}{{\mathbb{P}}}

\newcommand{\Z}{{\mathbb{Z}}}

\renewenvironment{proof}{\par \medskip \noindent
{\sc Proof:}}{}

\begingroup
\newtheorem{thm}{Theorem}   
\newtheorem{cor}[thm]{Corollary}     
\newtheorem{lemma}[thm]{Lemma}         
\newtheorem{prop}[thm]{Proposition}  

\newtheorem{defn}[thm]{Definition}   
   
\endgroup

\newenvironment{rem}[2]{\refstepcounter{thm} \label{#2} 
\par \medskip \noindent {\bf #1 \thethm .}}{\par \medskip}

\begin{document}


\pagenumbering{arabic}

\title[The Cohomology of Reflexive Sheaves]{The Cohomology of Reflexive Sheaves on Smooth Projective 3-folds}

\author[Peter Vermeire]{Peter Vermeire}

\address{Department of Mathematics, 214 Pearce, Central Michigan
University, Mount Pleasant MI 48859}

\email{verme1pj@cmich.edu}

\date{\today}

\begin{abstract}
We study the cohomology of reflexive rank $2$ sheaves on
smooth projective threefolds.  Applications are given to the moduli space of reflexive sheaves.
\end{abstract}

\maketitle

\section{Introduction}
We work over an algebraically closed field of characteristic $0$. 

There has been a tremendous amount of interest in recent years in the
study of curves on projective threefolds, and especially on the
general quintic in $\P^4$.  In this paper we continue the study, motivated by Hartshorne's work
\cite{hartvb,hart,hart2,hart3} on curves in
$\P^3$ via the Serre Correspondence, of reflexive sheaves on projective threefolds \cite{vermeirereflexive,vermeireparity,vermeiremoduli}.

Recall that a  coherent sheaf $\F$ is \textit{torsion-free} if the natural map of $\F$ to its double-dual $h:\F\rightarrow\F^{**}$ is injective, and that $\F$ is \textit{reflexive} if $h$ is an isomorphism. 
We refer the reader to \cite{hart} for basic properties of reflexive
sheaves.  Recall the following {\it Serre Correspondence} for
reflexive sheaves:
\begin{thm}\label{scorr}\cite[4.1]{hart}
Let $X$ be a smooth projective threefold, $M$ an invertible sheaf with $\H^1(X,M^*)=\H^2(X,M^*)=0$.  There is a one-to-one correspondence
between 
\begin{enumerate}
\item pairs $(\F,s)$ where $\F$ is a rank $2$ reflexive sheaf on $X$
  with $\wedge^2\F=M$ and $s\in\Gamma(\F)$ is a section whose zero
  scheme has codimension $2$
\item pairs $(Y,\xi)$ where $Y$ is a closed Cohen-Macaulay curve in
  $X$, generically a local complete intersection, and
  $\xi\in\Gamma(Y,\omega_Y^{\circ}\tensor\omega_X^*\tensor M^*)$ is a section
  which generates the sheaf $\omega_Y^{\circ}\tensor\omega_X^*\tensor M^*$
  except at finitely many points.
\end{enumerate}
Furthermore, $\c_3(\F)=2p_a(Y)-2+\c_1(X)\c_2(\F)-\c_1(\F)\c_2(\F)$.
\nopagebreak \hfill $\Box$ \par \medskip
\end{thm}

Note that if $\F$ is locally free, then the corresponding curve $Y$ is a local complete
intersection,  $\omega_Y^{\circ}\tensor\omega_X^*\tensor M^*\cong\O_Y$, $\xi$
is a non-zero section, and $c_3(\F)=0$.  In this case we say $Y$ is
{\em subcanonical}.

We also recall:
\begin{thm}\label{rr}
Let $\F$ be a coherent sheaf of rank $r$ on a smooth threefold $X$.  Then the Riemann-Roch formula is
\begin{eqnarray*}
\chi(X,\F) & = &
\frac{1}{6}\c_1(\F)^3-\frac{1}{2}\c_1(\F)\c_2(\F)-\frac{1}{2}\c_1(X)\c_2(\F)+\frac{1}{4}\c_1(X)\c_1(\F)^2
\\
 & & +\frac{1}{12}\c_1(X)^2\c_1(\F)+\frac{1}{12}\c_1(\F)\c_2(X)+\frac{r}{24}\c_1(X)\c_2(X)+\frac{1}{2}\c_3(\F) 
\end{eqnarray*}
\nopagebreak \hfill $\Box$ \par \medskip
\end{thm}

\begin{rem}{Example}{anticanonicaldet}
Let $X\subset\P^4$ be a smooth hypersurface of degree $r$ and let $\F$ be a rank $2$ reflexive sheaf on $X$ with $\wedge^2\F=\O_X(k)$.  Suppose $\F$ has a section whose zero set is a curve $C$ of degree $d$.  Then combining Proposition~\ref{rr} with the expression for $c_3(\F)$ in Theorem~\ref{scorr}, we obtain
$$\chi(X,\F)=\frac{1}{12}r(k+5-r)(2k^2+5k-kr+10-5r+r^2)+p_a(C)-1-kd$$
\nopagebreak \hfill $\Box$ \par \medskip
\end{rem}

\begin{prop}\label{localtoglobal}
Let $\F$ be a reflexive sheaf on a normal projective threefold $X$, $\G$ a sheaf of $\O_X$-modules.  Then there are isomorphisms

$\H^0(X,\sext^0_{\O_X}(\F,\G))=\Ext^0_{\O_X}(\F,\G)$

$\H^3(X,\sext^0_{\O_X}(\F,\G))=\Ext^3_{\O_X}(\F,\G)$

and an exact sequence
$$0\rightarrow \H^1(X,\sext^0_{\O_X}(\F,\G))\rightarrow \Ext^1_{\O_X}(\F,\G)\rightarrow \H^0(X,\sext^1_{\O_X}(\F,\G))$$
$$\rightarrow \H^2(X,\sext^0_{\O_X}(\F,\G))\rightarrow \Ext^2_{\O_X}(\F,\G)\rightarrow0$$
\end{prop}

\begin{proof}
This follows exactly as in \cite[2.5]{hart}; the key point is that a reflexive sheaf on a smooth threefold $X$ has homological dimension 1.
\nopagebreak \hfill $\Box$ \par \medskip
\end{proof}

In section 2 we give some elementary examples of the relationship between the structure and the cohomology of a reflexive sheaf $\F$.  In section 3 we give a simple example of how this structure is affected by the existence of a global section, as will exist under the Serre Correspondence above.  This motivates section 4 where we investigate the influence of global sections on the higher cohomology of $\F$.  Finally, in section 5 we give applications to the moduli space of reflexive sheaves as studied in \cite{vermeiremoduli}.  In particular, we have Theorems~\ref{fanotheorem} and ~\ref{cytheorem} which clarify results from \cite{vermeiremoduli} by simplifying the hypotheses.  Finally, we give the most general statement about the moduli of reflexive sheaves in Theorem~\ref{bigthm}.

\section{Basic Examples}

As a simple example of the connection between the cohomology and the structure of $\F$, we recall that for a locally free sheaf $\F$ and a very ample line bundle $L$, if $\h^i(X,\F\tensor L^n)=0$ for $i=1,2$, $n\in\Z$, we say $\F$ is \textit{$L$-ACM} (arithmetically Cohen-Macaulay), as in this case the associated curve $Y$ is arithmetically Cohen-Macaulay (in the embedding by $L$) if an only if $\F$ is $L$-ACM (\cite{cm1}).  If $\F$ is reflexive we have:

\begin{prop}\label{noreflexiveacm}
Let $\F$ be a rank $2$ reflexive sheaf on a smooth projective 3-fold $X$.  If $\H^2(X,\F\tensor L^n)=0$ for all $n\ll0$ and for some ample invertible sheaf $L$, then $\F$ is locally free.  
In particular, $\F$ reflexive and $L$-ACM implies that $\F$ is locally free.
\end{prop}

\begin{proof}
This is \cite[2.5.1]{hart}, where it is shown that $\H^2(X,\F\tensor L^n)=\c_3(\F)$ for all $n\ll0$, and that $\c_3(\F)=0$ if and only if $\F$ is locally free.
\nopagebreak \hfill $\Box$ \par \medskip
\end{proof}

As a second example, we look at a case where the Riemann-Roch formula becomes especially simple.

\begin{prop}\label{serreduality}
Let $\F$ be a rank $2$ reflexive sheaf on a smooth projective 3-fold $X$ with $\c_1(\omega_X)=\c_1(\F)$.  Then $\c_3(\F)=2\h^2(X,\F)-2\h^1(X,\F)$.  Hence the following are equivalent:
\begin{enumerate}
\item $\h^2(X,\F)\leq\h^1(X,\F)$
\item $\h^2(X,\F)=\h^1(X,\F)$
\item $\chi(X,\F)=0$
\item $\F$ is locally free
\end{enumerate}
\end{prop}

\begin{proof}
Substituting $\G=\O_X$ in Proposition~\ref{localtoglobal}, we have the exact sequence
$$0\rightarrow \H^1(X,\F^*\tensor\omega_X)\rightarrow \H^2(X,\F)^*\rightarrow \H^0(X,\sext^1_{\O_X}(\F,\omega_X))$$
$$\rightarrow \H^2(X,\F^*\tensor\omega_X)\rightarrow \H^1(X,\F)^*\rightarrow 0 $$
Further, by \cite[2.6]{hart} we have $\c_3(\F)=\h^0(X,\sext^1(\F,\omega_X))$.
Our hypotheses imply that $\F^*\tensor \omega_X=\F$ and the statements immediately follow.
\nopagebreak \hfill $\Box$ \par \medskip
\end{proof}

Lest the hypothesis in Proposition~\ref{serreduality} seem artificial, we point out that on a hypersurface $X$ every rank 2 reflexive sheaf can be twisted so that either $\c_1(\F(n))=\c_1(\omega_X)$ or $\c_1(\F(n))=\c_1(\omega_X(1))$.  In general, we have the following:

\begin{defn}
Let $X$ be a projective threefold.  A rank 2 reflexive sheaf $\F$ has \textbf{canonical parity} if $\wedge^2\F\tensor\omega_X^*=L^{\tensor 2}$ for some invertible sheaf $L$.
\nopagebreak \hfill $\Box$ \par \medskip
\end{defn}

\begin{rem}{Remark}{whycanonicalparity}
The point of the definition is that $\F$ has canonical parity if and only if there is an invertible sheaf $M$ such that $\c_1(\F\tensor M)=\c_1(\omega_X)$.  In fact, $M$ is unique and $M=\left(\wedge^2\F^*\tensor\omega_X\right)^{1/2}$.
\nopagebreak \hfill $\Box$ \par \medskip
\end{rem}

By Proposition~\ref{serreduality} we have
\begin{cor}
Let $\F$ be a rank $2$ reflexive sheaf with canonical parity on a smooth projective 3-fold $X$.  Then $\F$ is locally free if and only if 
$$\chi\left(X,\F\tensor\left(\wedge^2\F^*\tensor\omega_X\right)^{1/2}\right)=0  $$
\nopagebreak \hfill $\Box$ \par \medskip
\end{cor}

\section{Bounds on $\c_3$}

\begin{rem}{Notation}{S}  
For a fixed invertible sheaf $L$ and reflexive rank 2 sheaf $\F$ on a threefold $X$, we will use the notational simplification $S=\c_1(L)\c_2(\F)$.  Note that if $\F$ has a section whose zero scheme is a curve $C$ and if $L$ is very ample, then $S$ is simply the degree of $C$ in the embedding induced by $L$.
\nopagebreak \hfill $\Box$ \par \medskip
\end{rem}

In \cite{vermeirereflexive} our goal, for a fixed $X$, was to give a bound on
$\c_3(\F)$ in terms of $\c_1(\F)$ and $\c_2(\F)$.  Note that the formula
for $\c_3$ in Theorem~\ref{scorr} gives:
\begin{lemma}\label{sectionbound}\cite{vermeirereflexive}
Let $X$ be a smooth 3-fold, $L$ a very ample line bundle, $\F$ a
rank two reflexive sheaf with
$\h^1(\wedge^2\F^*)=\h^2(\wedge^2\F^*)=0$.  If $s\in\Gamma(\F)$ is a section whose zero
scheme is a curve, then $$\c_3(\F)\leq
S^2-3S+\c_1(X)\c_2(\F)-\c_1(\F)\c_2(\F)$$ 
\nopagebreak \hfill $\Box$ \par \medskip
\end{lemma}

A brute-force application of Riemann-Roch yields:

\begin{cor}\label{firstbound}
Let $L$ be a line bundle on a smooth projective 3-fold $X$, $\c_1(L)^3=d$, $\F$ a rank 2 coherent sheaf with $\c_1(\F)=\c_1(\O_X)$.  Suppose $n$ is such that $\H^2(X,\F\tensor L^n)=0$ and 
$$4n^3d+6n^2\c_1(X)\c_1(L)^2+2n\c_1(X)^2\c_1(L)+2n\c_2(X)\c_1(L)>$$ $$12nS+6\c_1(X)\c_2(F)-\c_1(X)\c_2(X)-6\c_3(F)$$
then $\F\tensor L^n$ has a section.  Suppose further that $\F$ is reflexive and $L$ very ample; then
$$\c_3(F)\leq (S-3-2n+2n^2d)(S+2n^2d)+\c_1(X)\c_2(F)+2n^2\c_1(X)\c_1(L)^2$$
\end{cor}

\begin{proof}
The first part follows directly from Theorem~\ref{rr}; the second part from Lemma~\ref{sectionbound}.
\nopagebreak \hfill $\Box$ \par \medskip
\end{proof}

The work cited above is primarily concerned with finding conditions under which $\H^2(X,\F\tensor L^n)=0$ so that one may apply Corollary~\ref{firstbound}.  For example:

\begin{rem}{Example}{oldbound}
In \cite[\S 3]{vermeirereflexive} it is shown that if $X\subset\P^4$ is a quintic hypersurface, and if $\F$ is a rank 2 semistable reflexive sheaf on $X$ with $\c_1(\F)=\c_1(\O_X)$, then we have $\H^2(X,\F(n))=0$ for 
$$\begin{cases}
n\geq 31 & \text{when } S\leq 19 \\
n\geq 4S-27+\frac{1}{2}\sqrt{60S-525} & \text{when } S\geq 20
\end{cases}$$
\nopagebreak \hfill $\Box$ \par \medskip
\end{rem}

Note that to a non-rational curve $C$ on a canonically trivial threefold $X$, we can associate a rank 2 reflexive sheaf $\F$ with $\wedge^2\F=\O_X$.  Because we know that $\F$ has a section, we don't need the first part of Corollary~\ref{firstbound}, and the second becomes:

\begin{prop}\label{foreshadow}
Let $X$ be a smooth projective threefold with $\H^1(X,\O_X)=\H^2(X,\O_X)=0$, and $\omega_X=\O_X$; let $C\subset X$ be a non-rational Cohen-Macaulay curve.  Then the sheaf $\F$ with $\wedge^2\F=\O_X$ associated to $C$ has $$2\chi(X,\F)=\c_3(\F)=2p_a(C)-2\leq S^2-3S$$
where $S$ is computed using any very ample invertible sheaf $L$.
\end{prop}

\begin{proof}
This follows immediately from Theorem~\ref{scorr}, Theorem~\ref{rr}, and Lemma~\ref{sectionbound}.
\nopagebreak \hfill $\Box$ \par \medskip
\end{proof}

Proposition~\ref{foreshadow} correctly suggests that the properties of reflexive sheaves associated to curves via the Serre Correspondence are not typical of the properties of reflexive sheaves in general.  While it is true that after twisting by a sufficient power of an ample line bundle any reflexive sheaf will have a section whose zero scheme is a curve, controlling this power can be quite delicate.  Furthermore, often our real interest is in studying the curves in a threefold rather than the reflexive sheaves the threefold carries.  Therefore, in the next section we study reflexive sheaves under the assumption that they admit a section.

\section{The influence of global sections on cohomology}

\begin{prop}\label{h0killsk2}
Suppose $X$ is a smooth threefold, $\F$ a rank $2$ reflexive sheaf.  If $\F$ has a section whose zero scheme is a curve $C$ and if $N$ is an invertible sheaf such that
\begin{enumerate}
\item $\H^2(X,\omega_X\tensor N)=0$
\item $\H^2(X,\wedge^2\F\tensor\omega_X\tensor N)=0$
\end{enumerate}
then $$\h^2(X,\F\tensor\omega_X\tensor N)\leq \h^1(C,\wedge^2\F\tensor\omega_X\tensor N\tensor\O_C)$$  Therefore, if further $\c_1(N)\c_2(\F)>\c_3(\F)$ then $\h^2(X,\F\tensor\omega_X\tensor N)=0$.

Similarly, if $\H^3(X,\omega_X\tensor N)=\H^3(X,\omega_X\tensor N\tensor\wedge^2\F)=0$, then $\H^3(X,\F\tensor\omega_X\tensor N)=0$.  
\end{prop}

\begin{proof}
The section of $\F$ induces the exact sequence
$$\ses{\omega_X\tensor N}{\F\tensor\omega_X\tensor N}{\I_C\tensor \wedge^2\F\tensor\omega_X\tensor N}$$
We have $\H^2(X,\omega_X\tensor N)=0$ hence 
$$\h^2(X,\F\tensor\omega_X\tensor N)\leq\h^2(X,\I_C\tensor \wedge^2\F\tensor\omega_X\tensor N)$$  From the standard ideal sheaf exact sequence, we have
$$\ses{\I_C\tensor \wedge^2\F\tensor\omega_X\tensor N}{\wedge^2\F\tensor\omega_X\tensor N}{\wedge^2\F\tensor\omega_X\tensor N\tensor\O_C}$$
The first part of the result now follows by the second hypothesis on $N$.
We know that $\omega_C^{\circ}=\wedge^2\F\tensor\omega_X\tensor\O_C(D)$ where $|D|=c_3(\F)$ (and in particular that if $c_3(\F)=0$ then $D=\emptyset$).  The hypothesis $\c_1(N)\c_2(\F)>\c_3(\F)$ gives $\h^1(C,\wedge^2\F\tensor\omega_X\tensor N\tensor\O_C)=0$ for degree reasons.

The last part of the Proposition follows similarly, but more easily.
\nopagebreak \hfill $\Box$ \par \medskip
\end{proof}

\begin{cor}\label{h2onhypersurface}
Let $\F$ be a reflexive sheaf of rank 2 on a smooth hypersurface $X\subset\P^4$.  If $\F$ has a section whose zero scheme is a curve $C$ and if $N$ is an invertible sheaf such that $\c_1(N)\c_2(\F)>\c_3(\F)$ (resp. $\geq\c_3(\F)$), then $\h^2(X,\F\tensor\omega_X\tensor N)=0$ (resp. $\leq 1$).  In particular, if $\F$ is locally free then $\H^2(X,\F(1)\tensor\omega_X)=0$ and $\h^2(X,\F\tensor\omega_X)=\h^0(C,\O_C)-1$.
\nopagebreak \hfill $\Box$ \par \medskip
\end{cor}

\begin{proof}
Except for the final statement, this follows from Proposition~\ref{h0killsk2}.  To get the equality $\h^2(X,\F\tensor\omega_X)=\h^0(C,\O_C)-1$, we note that $\H^2(X,\F\tensor\omega_X)=\H^1(X,\F^*)^*$ and so the sequence
$$\ses{\wedge^2\F^*}{\F^*}{\I_C}$$
completes the proof.
\nopagebreak \hfill $\Box$ \par \medskip
\end{proof}

We can verify the condition in Corollary~\ref{h2onhypersurface} as follows: Suppose $X$ has degree $r$ and let $N=\O_X(p)$ and $\wedge^2\F=\O_X(k)$; if $$p\geq\frac{2p_a(C)-2}{\deg(\O_C(1))}+5-r-k$$
then $\c_1(N)\c_2(\F)\geq\c_3(\F)$ where strict inequality holds in one if and only if it holds in the other.  To see this:
\begin{eqnarray*}
c_3(\F)&=&2p_a(C)-2+\c_2(\F)\c_1(X)-\c_2(\F)\c_1(\F)\\
&=&2p_a(C)-2-\c_2(\F)\c_1(\O_X(r-5+k)) \\
&=&2p_a(C)-2-(r-5+k)\deg(\O_C(1))\\
&\leq&p\deg(\O_C(1))\\
&=&\c_1(N)\c_2(\F)
\end{eqnarray*}

Using this we have:
\begin{cor}\label{justkillH2F}
Let $\F$ be a reflexive sheaf of rank 2 on a smooth hypersurface $X\subset\P^4$.  Suppose that $\F$ has a section whose zero scheme is a curve $C$ and that $\wedge^2\F=\O_X(k)$.  If $k\deg(C)>2p_a(C)-2$ then $H^2(X,\F)=0$.
\nopagebreak \hfill $\Box$ \par \medskip
\end{cor}

\begin{cor}\label{betterhypkillh0}
Let $\F$ be a reflexive sheaf of rank 2 on a smooth hypersurface $X\subset\P^4$ with $\wedge^2\F=\O_X(k)$.  If $\F$ has a section whose zero scheme is a curve $C$ and if $p>\max\{0,-k\}$
then $\H^3(X,\F(p)\tensor\omega_X)=\H^0(X,\F(-k-p))=0$.
\end{cor}

\begin{proof}
For convenience suppose $X\subset\P^4$ has degree $r$.  We use the proof of Proposition~\ref{h0killsk2} with $N=\O_X(p)$.  The section induces the exact sequence
$$\ses{\O_X(r-5+p)}{\F(r-5+p)}{\I_C(k+r-5+p)}$$
As long as $p>0$, we have $\H^3(X,\O_X(r-5+p))=0$ hence 
$$\H^3(X,\F(r-5+p))=\H^3(X,\I_C(k+r-5+p))$$  From the standard ideal sheaf exact sequence, we have
$$\ses{\I_C(k+r-5+p)}{\O_X(k+r-5+p)}{\O_C(k+r-5+p)}$$
and so
$$\H^3(X,\F(r-5+p))=\H^3(X,\O_X(k+r-5+p))$$
Finally, as long as $p>-k$ we have $\H^3(X,\O_X(k+r-5+p))=\H^0(X,\O(-k-p))=0$.  Therefore, $\H^3(X,\F(r-5+p))=\H^0(X,\F(-k-p))=0$.
\nopagebreak \hfill $\Box$ \par \medskip
\end{proof}

\begin{rem}{Example}{killsplitsheaves}
Let $C$ be a smooth non-rational curve on a smooth quintic hypersurface $X\subset\P^4$.  Then $C$ is the zero scheme of a section of a rank $2$ reflexive sheaf $\F$ with $\wedge^2\F=\O_X$.  From the sequence
$$\ses{\O_X}{\F}{\I_C}$$
we see $\H^0(X,\F(p))=0$ if and only if $p<0$.  Therefore, Corollary~\ref{betterhypkillh0} is the best possible in general.
\nopagebreak \hfill $\Box$ \par \medskip
\end{rem}

\begin{prop}\label{dualdies}
Let $\F$ be a rank $2$ reflexive sheaf on a smooth threefold $X\subset\P^4$ of degree $r$.  If $\wedge^2\F=\O_X(k)$ is ample and if $\F$ has a section whose zero scheme is a connected curve $C$ then $\H^0(X,\F^*)=\H^1(X,\F^*)=0$.  Further, we have an exact sequence
$$0\rightarrow\H^2(X,\F^*)\rightarrow \H^1(C,\O_C)\rightarrow \H^0(X,\O_X(r+k-5))^*$$
$$\rightarrow \H^3(X,\F^*)\rightarrow\H^0(X,\O_X(r-5))^*\rightarrow0  $$
hence
$$\chi(X,\F^*)=\h^2(X,\F^*)-\h^3(X,\F^*)=p_a(C)-\binom{r-1+k}{4}+\binom{k-1}{4}-\binom{r-1}{4}$$
\end{prop}

\begin{proof}
This follows immediately from the associated exact sequence $$\ses{\O_X(-k)}{\F^*}{\I_C}$$ along with some elementary calculation.
\nopagebreak \hfill $\Box$ \par \medskip
\end{proof}

\begin{rem}{Example}{ratcurveskillF*}
If $C\subset X$ is a rational curve that is not a line on a hypersurface of degree at most $3$, then taking $\wedge^2\F=\O_X(1)$ we see that $\H^2(X,\F^*)=\H^3(X,\F^*)=0$.
\nopagebreak \hfill $\Box$ \par \medskip
\end{rem}

\begin{cor}
Let $\F$ be a rank two reflexive sheaf on a smooth hypersurface $X\subset\P^4$ of degree $r\leq 4$.  Suppose that $\F$ has a section whose zero scheme is a connected curve $C$ and that $\wedge^2\F=\O_X(k)$ is effective.  If $H^0(X,\I_C(r-5+k))=0$ then 
$$h^2(X,\F^*)=\chi(X,\F^*)=p_a(C)-h^0(X,\O_X(r-5+k))$$
In particular, $p_a(C)\geq h^0(X,\O_X(r-5+k))$.
\nopagebreak \hfill $\Box$ \par \medskip
\end{cor}

\begin{proof}
This follows immediately from the fact that $\H^3(X,\F^*)=\H^0(X,\F\tensor\omega_X)^*$.
\nopagebreak \hfill $\Box$ \par \medskip
\end{proof}

\section{Moduli of sheaves}

In this section, we give an application to the moduli space of torsion-free sheaves using \cite{vermeiremoduli}.

\begin{prop}\label{F}
Let $X$ be a smooth projective threefold with $\H^2(X,\O_X)=0$, $\F$ a rank $2$ reflexive sheaf.  Suppose that $\F$ has a section whose zero locus is a curve $C$ and that $\wedge^2\F\tensor\O_C$ is non-special.  Then $\H^2(X,\F)=0$.
\end{prop}

\begin{proof}
We apply Proposition~\ref{h0killsk2} with $N=\omega_X^*$.  This yields $\h^2(X,\F)\leq\h^1(C,\wedge^2\F\tensor\O_C)=0$. 
\nopagebreak \hfill $\Box$ \par \medskip
\end{proof}

\begin{prop}\label{F*}
Let $X$ be a smooth projective threefold with $\H^2(X,\O_X)=0$, $\F$ a rank $2$ locally free sheaf.  Suppose that $\F$ has a section whose zero locus is a curve $C$ and that $\H^1(X,\I_C\tensor\wedge^2\F\tensor\omega_X)=0$.  Then $\H^2(X,\F^*)=0$.
\end{prop}

\begin{proof}
We have $\H^2(X,\F^*)=\H^1(X,\F\tensor\omega_X)^*$, and the result follows from the sequence $$\ses{\omega_X}{\F\tensor\omega_X}{\I_C\tensor\wedge^2\F\tensor\omega_X}$$
\nopagebreak \hfill $\Box$ \par \medskip
\end{proof}

\begin{thm}\label{fanotheorem}
Let $X$ be a smooth projective Fano threefold, $\F$ a stable rank $2$ locally free sheaf.  Suppose that $\F$ has a section whose zero scheme is a curve $C$ and suppose further that 
$\H^1(X,\I_C\tensor\wedge^2\F\tensor\omega_X)=0$.  If $\H^1(C,N_{C/X})=0$ then the (coarse) projective moduli space of semi-stable coherent rank $2$ torsion-free sheaves is smooth of dimension $$1-\frac{c_1(X)c_2(X)}{6}+\frac{c_1(X)\Delta(\F)}{2}$$
at the point corresponding to $\F$.
\end{thm}

\begin{proof}
Note that because $\F$ is locally free we have $\omega_C=\omega_X\tensor\wedge^2\F\tensor\O_C$.  In any case $\H^1(C,\wedge^2\F\tensor\O_C)=\H^0(C,\omega_X\tensor\O_C)^*=0$ because $\omega_X^*$ is ample, hence $\H^2(X,\F)=0$ by Proposition~\ref{F}.

The result now follows from Proposition~\ref{F*} and \cite{vermeiremoduli}.
\nopagebreak \hfill $\Box$ \par \medskip
\end{proof}

\begin{thm}\label{cytheorem}
Let $X$ be a smooth projective Calabi-Yau threefold with $\H^1(X,\O_X)=0$, $\F$ a stable rank $2$ locally free sheaf.  Suppose that $\F$ has a section whose zero scheme is a curve $C$ and suppose further that that $\H^1(X,\I_C\tensor\wedge^2\F)=0$.  If $\H^1(C,N_{C/X})=0$ then the (coarse) projective moduli space of semi-stable coherent rank $2$ torsion-free sheaves is smooth of dimension zero at the point corresponding to $\F$.
\end{thm}

\begin{proof}
In this case $\H^2(X,\F)=\H^1(X,\F^*)^*=0$, hence the result follows from Proposition~\ref{F*} and \cite{vermeiremoduli}.
\nopagebreak \hfill $\Box$ \par \medskip
\end{proof}

Finally, we are able, at a slight cost, to extend Theorem~\ref{fanotheorem} to the case where $\F$ is reflexive.  The main techincal result is:

\begin{prop}\label{fanothm2}
Let $X$ be a smooth projective threefold, $\F$ a rank $2$ reflexive sheaf.
Suppose that $\F$ has a section whose zero scheme is a curve $C$ and suppose further that 
\begin{enumerate}
\item $\H^2(X,\F)=0$
\item $\H^2(X,\F^*)=0$
\item $\H^0(C,\wedge^2\F\tensor\omega_X\tensor\O_C)=0$
\item $\H^0(C,\omega_X\tensor\O_C)=0$
\item $\H^1(C,N_{C/X})=0$
\end{enumerate}
Then $\Ext^2_{\O_X}(\F,\F)=0$.
\end{prop}

\begin{rem}{Remark}{betterthanbefore}
This was proved in \cite{vermeiremoduli} for $\F$ locally free assuming only conditions 1,2, and 5.  Note also that as we are primarily interested in Fano varieties, condition 4 is usually of no concern.  Note also that if $\wedge^2\F$ is effective then condition 3 implies condition 4.
\end{rem}

\begin{proof}
We let $\omega_X\tensor\wedge^2\F\tensor\O_C=\omega_C(-D)$ where $D$ is an effective divisor on $C$.

Applying $\Hom_{\O_X}(\cdot,\F)$ to the sequence
$$\ses{\O_X}{\F}{\I_C\tensor\wedge^2\F}$$
we have 
$$\cdots\rightarrow\Ext^2_{\O_X}(\I_C\tensor\wedge^2\F,\F)\rightarrow \Ext^2_{\O_X}(\F,\F)\rightarrow \Ext^2_{\O_X}(\O_X,\F)\rightarrow\cdots$$
where $\Ext^2_{\O_X}(\O_X,\F)=\H^2(X,\F)=0$ by hypothesis.  Therefore, it suffices to show that $\Ext^2_{\O_X}(\I_C\tensor\wedge^2\F,\F)=\Ext^2_{\O_X}(\I_C,\F^*)=0$.

Applying $\Hom_{\O_X}(\cdot,\F^*)$ to the sequence $$\ses{\I_C}{\O_X}{\O_C}$$ yields 
$$\cdots\rightarrow\Ext^2_{\O_X}(\O_X,\F^*)\rightarrow \Ext^2_{\O_X}(\I_C,\F^*)\rightarrow \Ext^3_{\O_X}(\O_C,\F^*)\rightarrow\cdots$$
where $\Ext^2_{\O_X}(\O_X,\F^*)=\H^2(X,\F^*)=0$ by hypothesis.  Therefore, it suffices to show that $\Ext^3_{\O_X}(\O_C,\F^*)=0$.

Applying $\Hom_{\O_X}(\O_C,\cdot)$ to the sequence $$\ses{\wedge^2\F^*}{\F^*}{\I_C}$$ we have $$\cdots\rightarrow\Ext^3_{\O_X}(\O_C,\wedge^2\F^*)\rightarrow \Ext^3_{\O_X}(\O_C,\F^*)\rightarrow \Ext^3_{\O_X}(\O_C,\I_C)\rightarrow0$$
where 
\begin{eqnarray*}
\Ext^3_{\O_X}(\O_C,\wedge^2\F^*)&=&\Ext^3_{\O_X}(\O_C\tensor\wedge^2\F\tensor\omega_X,\omega_X)\\
&=&\H^0(X,\O_C\tensor\wedge^2\F\tensor\omega_X)\\
&=&\H^0(C,\omega_C(-D))\\
&=&0
\end{eqnarray*}
hence is suffices to show $\Ext^3_{\O_X}(\O_C,\I_C)=0$.

Applying $\Hom_{\O_X}(\O_C,\cdot)$ to the sequence $$\ses{\I_C}{\O_X}{\O_C}$$ we have
$$\cdots\rightarrow\Ext^2_{\O_X}(\O_C,\O_C)\rightarrow \Ext^3_{\O_X}(\O_C,\I_C)\rightarrow \Ext^3_{\O_X}(\O_C,\O_X)\rightarrow\cdots$$
where $\Ext^3_{\O_X}(\O_C,\O_X)=\H^0(C,\omega_X\tensor\O_C)=0$ because $X$ is Fano.  Therefore it suffices to show that $\Ext^2_{\O_X}(\O_C,\O_C)=0$.

Applying $\Hom_{\O_X}(\cdot,\O_C)$ to the sequence $$\ses{\I_C}{\O_X}{\O_C}$$ we have
$$\cdots\rightarrow\Ext^1_{\O_X}(\I_C,\O_C)\rightarrow \Ext^2_{\O_X}(\O_C,\O_C)\rightarrow \Ext^2_{\O_X}(\O_X,\O_C)\rightarrow\cdots$$
where $\Ext^2_{\O_X}(\O_X,\O_C)=\H^2(X,\O_C)=0$ because $C$ is a curve.  Finally, it suffices to show that $\Ext^1_{\O_X}(\I_C,\O_C)=0$.  We do this by showing that $$\H^1(X,\shom_{\O_X}(\I_C,\O_C))=\H^0(X,\sext^1_{\O_X}(\I_C,\O_C))=0$$

First, note that
\begin{eqnarray*}
\shom_{\O_X}(\I_C,\O_C)&=&\shom_{\O_X}(\I_C,\shom_{\O_C}(\O_C,\O_C))\\
&=&\shom_{\O_C}(\I_C\tensor_{\O_X}\O_C,\O_C)\\
&=&\shom_{\O_C}(N^*_{C/X},\O_C)\\
&=&N_{C/X}
\end{eqnarray*}
and so $\H^1(X,\shom_{\O_X}(\I_C,\O_C))=0$ by hypothesis.

Applying $\shom_{\O_X}(\I_C,\cdot)$ to $$\ses{\I_C}{\O_X}{\O_C}$$ yields 
\begin{eqnarray*}
0&\rightarrow&\shom_{\O_X}(\I_C,\I_C)\rightarrow\shom_{\O_X}(\I_C,\O_X)\xrightarrow{0}\shom_{\O_X}(\I_C,\O_C)\\
&\rightarrow&\sext^1_{\O_X}(\I_C,\I_C)\rightarrow\sext^1_{\O_X}(\I_C,\O_X)\rightarrow\sext^1_{\O_X}(\I_C,\O_C)\\
&\rightarrow&\sext^2_{\O_X}(\I_C,\I_C)\rightarrow\cdots
\end{eqnarray*}
Now, $\sext^2_{\O_X}(\I_C,\I_C)=0$ because $\I_C$ has homological dimension $1$ since $C$ is Cohen-Macaulay.  We therefore have the sequence
$$0\rightarrow N_{C/X}\rightarrow\sext^1_{\O_X}(\I_C,\I_C)\rightarrow\sext^1_{\O_X}(\I_C,\O_X)\rightarrow\sext^1_{\O_X}(\I_C,\O_C)\rightarrow0$$
where all the sheaves are supported on $C$.  By \cite[p.137]{hart} we know further that $\sext^1_{\O_X}(\I_C,\O_X)=\omega_C\tensor\omega_X^*$.  By hypothesis we have $\H^1(C,\omega_C\tensor\omega_X^*)=\H^0(C,\omega_X\tensor\O_C)=0$, hence $\H^1(X,\sext^1_{\O_X}(\I_C,\O_C))=0$ because $C$ has dimension $1$.
\nopagebreak \hfill $\Box$ \par \medskip
\end{proof}

Combining Proposition~\ref{fanothm2} with the results of \cite{vermeiremoduli}, we have:

\begin{thm}\label{bigthm}
Let $\F$ be a stable reflexive sheaf of rank $2$ on a smooth projective threefold $X\subset\P^n$ and assume that either $\omega_X^*$ is effective or that there exists an $n\in\Z$ such that $\H^0(X,\F\tensor\omega_X^n)\neq0$ and $\H^0(X,\F\tensor\omega_X^{n+1})=0$. 
Suppose further that $\F$ has a section whose zero scheme is a curve $C$, and that
\begin{enumerate}
\item $\H^2(X,\F^*)=0$
\item $\H^2(X,\F)=0$
\item $\H^0(C,\wedge^2\F\tensor\omega_X\tensor\O_C)=0$
\item $\H^0(C,\omega_X\tensor\O_C)=0$
\item $\H^1(C,N_{C/X})=0$
\end{enumerate}
Then the (coarse) projective moduli space of semi-stable coherent rank $2$ torsion-free sheaves at the point corresponding to $\F$ is smooth of dimension $$1-\frac{c_1(X)c_2(X)}{6}+\frac{c_1(X)\Delta(\F)}{2}$$
where $\Delta(\F)=4\c_2(\F)-\c_1^2(\F)$.
\nopagebreak \hfill $\Box$ \par \medskip
\end{thm}

A particularly nice example is:

\begin{cor}\label{fanocor}
Let $X$ be a smooth projective Fano threefold, $\F$ a stable rank $2$ reflexive sheaf with $\wedge^2\F$ big and nef.  Suppose that $\F$ has a section whose zero scheme is a rational curve $C$. If $\H^1(C,N_{C/X})=0$ then the (coarse) projective moduli space of semi-stable coherent rank $2$ torsion-free sheaves is smooth of dimension $$1-\frac{c_1(X)c_2(X)}{6}+\frac{c_1(X)\Delta(\F)}{2}$$ at the point corresponding to $\F$.
\end{cor}

\begin{rem}{Remark}{altforvanishing}
It is immediate that the condition $\wedge^2\F$ nef implies $\H^1(C,\wedge^2\F\tensor\O_C)=0$, which is not hard to see is equivalent to $2+\deg(N_{C/X})>c_3(\F)$. Further, if $\picard(X)=\Z$ then stability of $\F$ implies $\wedge^2\F$ is big and nef.
\end{rem}

\begin{proof}
We let $\omega_X\tensor\wedge^2\F=\omega_C(-D)$ where $D$ is an effective divisor on $C$.
By hypothesis $\H^1(C,\wedge^2\F\tensor\O_C)=0$, hence $\H^2(X,\F)=0$ by Proposition~\ref{F}.

From the sequence $$\ses{\wedge^2\F^*}{\F^*}{\I_C}$$ we have $\H^2(X,\wedge^2\F^*)=0$ by hypothesis and $\H^2(X,\I_C)=0$ because $C$ is rational.  Hence $\H^2(X,\F^*)=0$.

The result follows from Theorem~\ref{bigthm}.
\nopagebreak \hfill $\Box$ \par \medskip
\end{proof}

\end{document}